\documentclass[a4paper,12pt]{amsart}
\usepackage[
    hmargin=3cm,
    headheight=15pt
]{geometry}
\usepackage{amsmath,amssymb,amsthm,mathtools}
\usepackage{dsfont}
\usepackage{bbm}
\usepackage[nodate]{datetime}
\usepackage{color}
\usepackage{subfig}
\usepackage{fancyhdr}
\usepackage{enumerate}
\usepackage{todonotes}

\pdfoutput=1
 	
\newtheorem{theorem}{Theorem}[section]
\newtheorem*{theorem*}{Theorem}
\newtheorem{lemma}{Lemma}[section]

\newtheorem{definition}{Definition}[section]

\theoremstyle{remark}
\newtheorem{example}{Example}[section]
\newtheorem{remark}{Remark}[section]

\title{On certain correlations into the divisor problem}
\author{Alexandre Dieguez}
\providecommand{\keywords}[1]{\textbf{\textit{Index terms---}} #1}
\address{Departamento de Matem\'atica, Universidade Federal de Minas Gerais, Av. Ant\^onio Carlos, 6627, CEP 31270-901, Belo Horizonte, MG, Brazil.}
\email{alexdieguez33@gmail.com}
\date{\today}
\subjclass[2010]{}
\begin{document}
\begin{abstract}
For a fixed irrational \(\theta > 0\) with a prescribed irrationality measure function, we study the correlation \(\int_1^X \Delta(x) \Delta(\theta x)  dx\), where \(\Delta\) is the Dirichlet error term in the divisor problem. When \(\theta\) has a finite irrationality measure, it is known that decorrelation occurs at a rate expressible in terms of this measure. Strong decorrelation occurs for all positive irrationals, except possibly Liouville numbers. We show that for irrationals with a prescribed irrationality measure function \(\psi\), decorrelation can be quantified in terms of \(\psi^{-1}\).
\end{abstract}
\maketitle
\section{Introduction}
\subsection{The Dirichlet Divisor Problem}
One of the oldest open problems in analytic number theory, dating back to Dirichlet in the 19th century, is the \textit{Dirichlet divisor problem}. Applying Dirichlet's hyperbola method yields the asymptotic formula for the summatory function of the divisor function \(\tau\):
\begin{equation}\label{equation hyperbola}
\sum_{n \leq x} \tau(n) = x \log x + (2\gamma - 1) x + \Delta(x),
\end{equation}
where \(\gamma \approx 0.57721\) is the Euler-Mascheroni constant and \(\Delta\) is the error term. The problem seeks the least \(\alpha>0\) such that \(\Delta(x)=O(x^{\alpha+\varepsilon})\) is satisfied for every \(\varepsilon>0\). Dirichlet \cite{Dirichlet_2012} used his hyperbola method to show that \(\alpha\leq 1/2\).

During the past century, considerable effort has been devoted to improving this bound, both through sharper estimates for \(\alpha\) and improved bounds for \(\Delta\). Vorono\"{\i} \cite{Voronoi_1904a,Voronoi_1904b} established \(\Delta(x)=O(x^{1/3}\log x)\) in 1904, giving \(\alpha \leq 1/3\). A breakthrough occurred in 1916 when Hardy \cite{Hardy_1916}, building on Landau's unpublished work, proved \(\Delta(x)=\Omega_{\pm}(x^{1/4})\), establishing \(1/4\) as the first lower bound for \(\alpha\) and leading to the conjecture that \(\alpha=1/4\).

Throughout the twentieth century, incremental improvements to the upper bound for \(\alpha\) were made by van der Corput, Littlewood, Walfisz, Tsung-tao, Hua, Kolesnik, Vinogradov, Iwaniec, and Mozzochi. The current record (Huxley, 2003 \cite{Huxley_2003}) is \(\alpha \leq 131/416\).

Although most of the work focused on upper bounds, gradual progress has also been made on the \(\Omega_{-}\) and \(\Omega_{+}\) bounds for \(\Delta\) since Hardy's result. C\'orradi and K\'atai \cite{Corradi_Katai_1967} improved the \(\Omega_{-}\) bound in 1967, while Hafner \cite{Hafner_1981} improved the \(\Omega_{+}\) bound in 1981. Soundararajan \cite{Soundararajan_2003} achieved a major refinement in 2003:
\[
\Delta(x) = \Omega \left( (x \log x)^{1/4} \frac{(\log \log x)^{(3/4)(2^{4/3} - 1)}}{(\log \log \log x)^{5/8}} \right).
\]
This result, together with Huxley's upper bound, constitutes a key milestone in 21st-century progress on the problem.

\subsection{Autocorrelation in the Divisor Problem}
In \cite{Aymone_etal_2024}, Aymone, Maiti, Ramar\'e, and Srivastav investigate
\begin{equation}\label{equation mainintegral}
    \int_1^X \Delta(x) \Delta(\theta x)  dx, \quad \theta > 0.
\end{equation}
The asymptotic behavior depends on whether \(\theta\) is rational, 
irrational with a finite irrationality measure, or a Liouville number. 
When \(\theta = a/b\) is rational, they establish the non-vanishing correlation:
\[
\lim_{X \to +\infty} \frac{1}{X^{3/2}} \int_1^X \Delta(x) \Delta\left(\frac{ax}{b}\right)  dx = \frac{C}{\sqrt{\gcd(a, b)}} \varphi\left( \frac{\text{lcm}(a, b)}{\gcd(a, b)} \right),
\]
where \(C>0\) is a constant and \(\varphi\) denotes a positive multiplicative function. For irrational \(\theta\) with finite irrationality measure \(\eta+1\) – which means that for some \(C > 0\) (typically \(C = 1\)) and all fixed \(\delta > 0\), the inequality
\[
|n - m \theta| \geq \frac{C}{m^{\eta + \delta}}
\]
is satisfied for all but finitely many rationals \(n/m\) – they prove
\[
\int_1^X \Delta(x) \Delta(\theta x)  dx \ll_{\varepsilon} X^{3/2 - 1/(18\eta) + \varepsilon}.
\]
A similar decorrelation result was independently obtained by Ivi\'c and Zhai \cite{Ivic_Zhai_2020}.

For Liouville numbers (irrationals that admit no finite irrationality measure and are well approximated by rationals), the integral satisfies \(o\left(X^{3/2}\right)\). Since Khinchin \cite{Khinchin_1924} established that \textit{almost all} irrationals have irrationality measure 2, the integral in \eqref{equation mainintegral} strongly decorrelates at the rate \(X^{3/2-1/18+\varepsilon}\) for almost all irrational \(\theta > 0\), with the Liouville numbers being exceptions.

\subsection{Liouville Numbers, \(\psi\)-Approximable Numbers, and Main Result}
Recall that a \textit{Liouville number} is an irrational \(\theta\) without a finite irrationality measure, which means that for every \(\eta>0\), there exist infinitely many rationals \(n/m\) satisfying
\[
|n - m\theta| < \frac{1}{m^{\eta}}.
\]

This work focuses on improving the asymptotic estimate for \eqref{equation mainintegral} when \(\theta > 0\) belongs to specific classes of Liouville numbers characterized by a prescribed increasing function \(\psi\), as defined below.

Building on Jarn\'{\i}k's foundational work \cite{Jarnik_1931} and Sondow's extensions \cite{Sondow_2003}, we generalize the concept of irrationality measure to an \textit{irrationality measure function}. Jarn\'{\i}k demonstrated that irrationals with a prescribed measure function \(\psi\) can be constructed under appropriate conditions, while Sondow introduced the \textit{irrationality base} to quantify rational approximation quality for Liouville numbers.

We note that the literature lacks a standard definition of irrationality measure function. Some sources require \(\psi\) to be a bivariate function increasing in both arguments, while Jarn\'{\i}k \cite{Jarnik_1931} uses a univariate function. We adopt the latter approach, in accordance with Jarn\'{\i}k's Theorem (see Theorem \ref{theorem Jarnik} below) and Sondow's irrationality base concept.

For an irrational \(\theta\) and an integer \(m\geq 1\), define \(||m\theta|| := \min_{n \in \mathbb{Z}} |m \theta - n|\) (the distance to the nearest integer). We say \(\theta\) has \textit{irrationality measure function} \(\psi\) if it is not \textit{approximable to order \(\psi\)}:

\begin{definition}
    Let \(\psi: [1, +\infty) \to (0, +\infty)\) be increasing. An irrational \(\theta\) is \textit{approximable to order \(\psi\)} (\(\psi\)-approximable) if 
    \[ 
    ||m\theta|| < \frac{1}{\psi(m)}
    \]
    for infinitely many integers \(m \geq1 \). Otherwise, \(\theta\) is not approximable to order \(\psi\).
\end{definition}

Khinchin's Theorem \cite{Khinchin_1924} states that if \(x/\psi(x)\) decreases and 
\[
\sum_{m\geq 1} \frac{1}{\psi(m)} < +\infty,
\] 
then almost all irrationals are not approximable to order \(\psi\) (i.e., they have the irrationality measure function \(\psi\)). Beyond this almost-all case, our main objective is to improve the upper bound for \eqref{equation mainintegral} for individual irrationals \(\theta > 0\) not approximable to order \(\psi\), which includes certain Liouville numbers.

With these definitions, we state our main result. Its scope includes Liouville numbers whose existence is guaranteed by Jarn\'{\i}k's Theorem (see Preliminaries section).

\begin{theorem}\label{theorem main}
    Let \(\theta > 0\) be irrational and \(\psi: [1, +\infty) \to (0, +\infty)\) increasing with \(\psi^{-1}(x) = O(x^{1/4})\). Suppose that there exists \(C > 0\) such that 
    \[
    ||m\theta|| \geq \frac{C}{\psi(m)}
    \] 
    for all sufficiently large integers \(m \geq 1\). Then
\[
\int_1^X \Delta(x) \Delta(\theta x)  dx \ll_{\varepsilon} \frac{X^{3/2}}{\psi^{-1}(X^{1/4})^{3/2 - \varepsilon}}.
\]
\end{theorem}

Since Khinchin's Theorem \cite{Khinchin_1924} states that almost all irrationals have irrationality measure 2, we provide an improved estimate covering the generic case:

\begin{theorem}\label{theorem second}
    For an irrational \(\theta > 0\) with finite irrationality measure \(\eta + 1 < 5/2\),
\[
\int_1^X \Delta(x) \Delta(\theta x)  dx \ll_{\varepsilon} X^{3/2 - 1/8 + \varepsilon}.
\]
\end{theorem}

While Theorem \ref{theorem main} addresses general \(\psi\), imposing \(\frac{1}{\psi(x)} = o(x^{-1})\) ensures via Jarn\'{\i}k's Theorem the existence of an irrational \(\theta > 0\) approximable to order \(\psi\), but not to order \(C\psi\) for any \(C > 1\).

\begin{example}\label{example 1.1}
    Sondow \cite{Sondow_2003} defines the \textit{irrationality base} of an irrational \(\theta\) as the least \(\beta \in [1, +\infty)\) satisfying: for every \(\varepsilon > 0\),
\[
||m\theta|| \geq \frac{1}{(\beta + \varepsilon)^m}
\] 
for all sufficiently large integers \(m \geq 1\). By Proposition 1 \cite{Sondow_2003}, \(\tau_\beta = \theta\) with finite irrationality base \(\beta > 1\) is a Liouville number, while nonexistent \(\beta\) (\(\beta=+\infty\)) defines a \textit{super Liouville number}.

In our framework, such \(\tau_\beta\) is approximable to order \(\beta^x\), but not to order \((\beta + \varepsilon)^x\) for any \(\varepsilon > 0\). Theorem 4 in \cite{Sondow_2003} constructs explicit Liouville numbers: for rational \(\beta = a/b > 1\) in the lowest terms, let
\[
\tau_\beta = \frac{1}{\beta} + \frac{1}{\beta^a} + \frac{1}{\beta^{a^a}} + \cdots,
\]
which is not approximable to order \(\beta^x\). Applying Theorem \ref{theorem main} with \(\psi(x) = (\beta + 1)^x\) and noting \(\psi^{-1}(x) = \log_{\beta+1} x = \frac{\log x}{\log(\beta+1)}\), we obtain for \(\tau_\beta\):
\[
\int_1^X \Delta(x) \Delta(\tau_\beta x)  dx \ll_{\varepsilon} \frac{X^{3/2}}{(\log X)^{3/2 - \varepsilon}}.
\]
This bound is satisfied for all \(\tau_\beta\) with finite irrationality base \(\beta \geq 1\), since \(\tau_\beta\) is not approximable to order \((\beta + 1)^x\).
\end{example}

\begin{example}\label{example 1.2}
    Define \(\psi(x) = \exp(\exp{x})\). By Jarn\'{\i}k's Theorem, there exists an irrational \(\theta > 0\) approximable to order \(\psi\), but not to order \(2\psi\). This qualifies \(\theta\) as a super Liouville number (Example \ref{example 1.1}), since for any \(\beta > 1\), 
    \[
    \frac{1}{\psi(x)} = o(\beta^{-x}) \implies ||m\theta|| < \frac{1}{\beta^m}
    \]
    for infinitely many integers \(m \geq 1\). Theorem \ref{theorem main} with \(\psi^{-1}(x) = \log \log x\) yields:
\[
\int_1^X \Delta(x) \Delta(\theta x)  dx \ll_{\varepsilon} \frac{X^{3/2}}{(\log \log X)^{3/2 - \varepsilon}}.
\]
\end{example}

\subsection{Proof Strategy}
To evaluate \eqref{equation mainintegral}, we employ a refined closed-form expression for \(\Delta\) originally due to Vorono\"{\i} \cite{Voronoi_1904a,Voronoi_1904b} and improved by Lau and Tsang \cite{Lau_Tsang_1995}. For any \(N>0\) and \(\varepsilon>0\),
\begin{equation}\label{equation voronoideltasum}
    \Delta(x) = \frac{x^{1/4}}{\sqrt{2}\pi} \sum_{n \leq N} \frac{\tau(n)}{n^{3/4}} \cos\left(4\pi \sqrt{nx} - \pi/4 \right) + R_N(x),
\end{equation}
with remainder term \( R_N(x) = O\left(x^\varepsilon + x^{1/2 + \varepsilon} N^{-1/2}\right)\).

After expanding the integral, we apply the Cauchy-Schwarz inequality combined with Tong's mean square estimate for \(\Delta\) \cite{Tong_1956} to bound the resulting double sum, following \cite{Aymone_etal_2024}. Through careful selection of parameters \(N = N(X)\) and \(T = T(X)\), we divide the sum over the indices \(m,n\) using a diagonal cutoff at \(T\).

Complementing the approach in \cite{Aymone_etal_2024}, our proof's key innovation lies in exploiting the subtle Diophantine properties of \(\theta\), particularly through Legendre's criterion. This fundamental result provides a necessary condition for rational approximations to be convergents of \(\theta\).

\section{Background on Continued Fractions}
All the following basic results on continued fractions can be found in either Khinchin's book \cite{Khinchin1997} or Tenenbaum's \cite{Tenenbaum_2015}.  
\subsection{Expansion and Convergents}
Let \(\theta\) be an irrational number. Its \textit{simple} (or \textit{regular}) \textit{continued fraction expansion} is:
\[
\theta = a_0 + \cfrac{1}{a_1 + \cfrac{1}{a_2 + \cfrac{1}{a_3 + \cdots}}} =: [a_0; a_1, a_2, a_3, \dots],
\]
where \(a_0 = \lfloor \theta \rfloor\) and \(a_k\) is a positive integer for \(k \geq 1\). The sequence \((a_k)_{k=0}^\infty\) is generated recursively via:
\[
\alpha_0 := \theta, \quad a_k := \lfloor \alpha_k \rfloor, \quad \alpha_{k+1} := \frac{1}{\alpha_k - a_k} \quad (k \geq 0).
\]

\textbf{Convergents.} The \textit{\(k\)}th \textit{convergent} of \(\theta\) is the rational number:
\[
\frac{n_k}{m_k} := [a_0; a_1, \dots, a_k]
\]
obtained by truncating the expansion after \(k\) terms. The numerators \(n_k\) and the denominators \(m_k\) satisfy the recurrence relations:
\begin{align*}
    n_k &= a_k n_{k-1} + n_{k-2}, \\
    m_k &= a_k m_{k-1} + m_{k-2}
\end{align*}
for \(k \geq 2\), with initial conditions:
\[
n_0 = a_0, \quad n_1 = a_0 a_1 + 1, \quad m_0 = 1, \quad m_1 = a_1.
\]

\subsection{Key Properties of Convergents}
The convergents \(\left( n_k/m_k \right)_{k=0}^\infty\) provide optimal rational approximations to \(\theta\):

\textbf{Best Approximation.} For any \(k \geq 1\),
\[
\frac{1}{m_{k+1}m_k+m_k^2} < \left| \theta - \frac{n_k}{m_k} \right| < \frac{1}{m_{k+1}m_k}.
\]
The quantity \( |m_k \theta - n_k| \) minimizes the approximation error for all fractions with denominator \(\leq m_k\):
\[
||m_k\theta||=|m_k \theta - n_k| = \min_{\substack{n \in \mathbb{Z} \\ 1 \leq m \leq m_k}} |m \theta - n|.
\]
In fact, it is minimal even among the denominators \(< m_{k+1}\):
\[
|m_k \theta - n_k| < |m \theta - n| \quad \text{for all} \quad 1 \leq m < m_{k+1}, \quad m \neq m_k.
\]
    
\textbf{Alternating sign and Determinant Identity.} The approximation alternates in sign:
\[
\frac{n_{2k-2}}{m_{2k-2}} < \theta < \frac{n_{2k-1}}{m_{2k-1}}, \quad \frac{n_{k}}{m_{k}}-\frac{n_{k-1}}{m_{k-1}} = \frac{(-1)^{k-1}}{m_{k}m_{k-1}} \quad (k \geq 1).
\]
From the last equality:
\[
n_k m_{k-1} - n_{k-1} m_k = (-1)^{k-1} \quad (k \geq 1),
\]
implying \(\gcd(m_k, n_k) = 1\).

\textbf{Denominator Growth.} The denominators satisfy \(m_k \geq F_{k+1}\) (\textit{Fibonacci} sequence) for all \(k \geq 0\), which implies exponential growth (see Remark \ref{remark 4.1} below).

\section{Preliminaries}
We state a lemma, originally due to Legendre, that establishes a necessary condition for a rational number to be a convergent for an irrational number \(\theta\). We then state additional lemmas from \cite{Aymone_etal_2024} concerning estimates for an integral and a function that will play a role in our proof. Finally, we state Jarn\'{\i}k's Theorem.
\begin{lemma}\label{lemma Legendre}(Legendre's Criterion, \cite{Legendre_2009})
    Let \(\theta\) be irrational and let \(n/m\) rational satisfy
\[
|n - m\theta| < \frac{1}{2m},
\]
    Then \(n/m\) is a convergent of the continued fraction expansion of \(\theta\).
\end{lemma}

By adopting the notation used in \cite{Aymone_etal_2024}, we state the following lemmas:
\begin{lemma}\label{lemma integral}(as Lemma 4.1 in \cite{Aymone_etal_2024})
    Let \(a>0\). Then
\[
\int x^{2} \cos(ax)  dx = x^{2} \frac{\sin(ax)}{a} + 2x \frac{\cos(ax)}{a^2} - 2 \frac{\sin(ax)}{a^3} + C,
\]
    where \(C\) is a constant. Furthermore, for \(X\geq 1\),
\[
\int\limits_{1}^{X} x^{2} \sin(ax)  dx \ll \frac{X^{2}}{a}, \quad \int\limits_{1}^{X} x^{2} \cos(ax)  dx \ll \frac{X^{2}}{a}.
\]
\end{lemma}
\begin{lemma}\label{lemma functionLambda}(part of Proposition 1.1 in \cite{Aymone_etal_2024})
    The function \(\Lambda:\mathbb{R}\to\mathbb{R}\) is defined as:
    \[
\Lambda(x) =
\begin{cases} 
\frac{1}{3}, & \text{if } x = 0, \\
\frac{\sin(x)}{x} + 2 \frac{\cos(x)}{x^2} - 2 \frac{\sin(x)}{x^3}, & \text{if } x \neq 0.
\end{cases}
\]
It is continuous and bounded. Moreover, for \(X\geq 1\) and fixed \(a>0\),
\[
\Lambda\left(a\sqrt{X}\right) = \frac{1}{X^{3/2}} \int_1^{\sqrt{X}} x^2 \cos(ax)  dx + \frac{\Lambda(a)}{X^{3/2}}.
\]
\end{lemma}
\begin{theorem}\label{theorem Jarnik}(Jarn\'{\i}k; Satz 6 in \cite{Jarnik_1931})
    Let \(\psi:[1, +\infty)\to(0, +\infty)\) be increasing with \(\frac{1}{\psi(x)} = o(x^{-1})\). Then there exists an irrational \(\theta > 0\) that is approximable to order \(\psi\), but not to order \(C\psi\) for any \(C>1\).
\end{theorem}

\section{Proof of The Main Result}
We now begin to prove our main result as in Theorem \ref{theorem main}:
\begin{proof}[Proof of Theorem 1.1.]
    Let \(\varepsilon > 0 \), which will always be sufficiently small and simplified in the upcoming expressions, and \(N>0\), which will be selected later. For \(1\leq x\leq X\), we have, by Vorono\"{\i}'s formula for \(\Delta\) in \eqref{equation voronoideltasum}:
\[
\Delta(x) = Q_N(x) + R_N(x),
\]
where
\[
Q_N(x) = \frac{x^{1/4}}{\sqrt{2}\pi} \sum_{n \leq N} \frac{\tau(n)}{n^{3/4}} \cos\left(4\pi \sqrt{nx} - \pi/4\right),
\]
and
\[
R_N(x) = O\left( x^\varepsilon + x^{1/2 + \varepsilon}N^{-1/2} \right).
\]
Now, we expand the integral in \eqref{equation mainintegral}:
\begin{align*}
\int_1^X \Delta(x) \Delta(\theta x)  dx 
&= \int_1^X Q_N(x) Q_N(\theta x)  dx 
+ \int_1^X \Delta(x) R_N(\theta x)  dx \\
&\quad + \int_1^X R_N(x) \Delta(\theta x)  dx 
- \int_1^X R_N(x) R_N(\theta x)  dx \\
&= \int_1^X Q_N(x) Q_N(\theta x)  dx \\
&+ O\left( X^{5/4+\varepsilon} + X^{7/4 + \varepsilon}N^{-1/2} + X^{2 + \varepsilon}N^{-1} \right).
\end{align*}
where the last three integrals of the first equality are estimated using the Cauchy-Schwarz inequality and \(\int_1^X \Delta(x)^2  dx \ll X^{3/2}\) (see Tong's result in \cite{Tong_1956}).
Next, we perform the change of variable \(u=x^{1/2}\), resulting in the following:
\begin{align*}
I_{\theta}(X) 
&= \int_1^X \Delta(x) \Delta(\theta x)  dx \\
&= \frac{1}{\pi^2} \sum_{m, n \leq N} \frac{\tau(m) \tau(n)}{(mn)^{3/4}} \int_1^{\sqrt{X}} u^2 \cos\left( 4\pi \sqrt{n} u - \pi/4 \right) \cos\left( 4\pi \sqrt{m\theta} u - \pi/4 \right)  du \\
&\quad + O\left( X^{5/4+\varepsilon} + X^{7/4 + \varepsilon}N^{-1/2} + X^{2 + \varepsilon}N^{-1} \right).
\end{align*}
We now apply the sum-to-product formula \(2\cos(u)\cos(v)=\sin(u+v)+\cos(u-v)\) to express the integral as:
\begin{align*}
&\frac{1}{2\pi^2} \sum_{m, n \leq N} \frac{\tau(m) \tau(n)}{(mn)^{3/4}} \int_1^{\sqrt{X}} u^2 \sin\left( 4\pi \left( \sqrt{m\theta} + \sqrt{n} \right) u \right)  du \\
&\quad + \frac{1}{2\pi^2} \sum_{m, n \leq N} \frac{\tau(m) \tau(n)}{(mn)^{3/4}} \int_1^{\sqrt{X}} u^2 \cos\left( 4\pi \left( \sqrt{m\theta} - \sqrt{n} \right) u \right)  du \\
&\quad + O\left( X^{5/4+\varepsilon} + X^{7/4 + \varepsilon}N^{-1/2} + X^{2 + \varepsilon}N^{-1} \right).
\end{align*}
Next, we invoke Lemma \ref{lemma integral} and use \(\tau(m) = o(m^\varepsilon) = o(N^\varepsilon)\) (see Montgomery and Vaughan \cite[pp.~55--56]{Montgomery_Vaughan_2006}):
\begin{align*}
&\frac{1}{2\pi^2} \sum_{m, n \leq N} \frac{\tau(m) \tau(n)}{(mn)^{3/4}} \int_1^{\sqrt{X}} u^2 \sin\left( 4\pi \left( \sqrt{m\theta} + \sqrt{n} \right) u \right) du \\
&\quad \ll X N^\varepsilon \sum_{m, n \leq N} \frac{1}{(mn)^{3/4} \left( \sqrt{m\theta} + \sqrt{n} \right)} \\
&\quad \ll X N^\varepsilon \sum_{m, n \leq N} \frac{1}{m^{5/4} n^{3/4}} \\
&\quad \ll X N^{\varepsilon}.
\end{align*}
Then, with Lemmas \ref{lemma integral} and \ref{lemma functionLambda}, by defining \( a_{m,n} := 4\pi \left( \sqrt{m\theta} - \sqrt{n} \right) \), we simplify the asymptotic expression for \( I_{\theta}(X) \) as follows:
\begin{align*}
I_{\theta}(X) &= \frac{1}{2\pi^{2}} \sum_{m, n \leq N} \frac{\tau(m)\tau(n)}{(mn)^{3/4}} \int_1^{\sqrt{X}} u^2 \cos \left( 4\pi \left( \sqrt{m\theta} - \sqrt{n} \right) u \right)  du \\
&\quad + O\left( X^{5/4+\varepsilon} + X^{7/4 + \varepsilon}N^{-1/2} + X^{2 + \varepsilon}N^{-1} + X N^{\varepsilon} \right) \\
&= \frac{X^{3/2}}{2\pi^{2}} \sum_{m, n \leq N} \frac{\tau(m)\tau(n)}{(mn)^{3/4}} \Lambda \left( a_{m,n} \sqrt{X} \right) \\
&\quad + O\left( X^{5/4+\varepsilon} + X^{7/4 + \varepsilon}N^{-1/2} + X^{2 + \varepsilon}N^{-1} + X N^{\varepsilon} \right).
\end{align*}
If we define:
    \[
J_{\theta}(X) := \frac{X^{3/2}}{2\pi^2} \sum_{m, n \leq N} \frac{\tau(m) \tau(n)}{(mn)^{3/4}} \Lambda \left( a_{m,n} \sqrt{X} \right),
\]
we achieve a final expression for \(I_{\theta}(X)\) as follows:
\[
    I_{\theta}(X) = J_{\theta}(X) + \tilde{R}_N(X),
\]
where 
\[
    \tilde{R}_N(X) = O \left( X^{5/4+\varepsilon} + X^{7/4 + \varepsilon}N^{-1/2} + X^{2 + \varepsilon}N^{-1} + X N^{\varepsilon} \right).
\]
To estimate \(I_{\theta}(X)\), we set the parameters \(N := X^{3/4}\) and 
\[
T := \frac{\pi}{\sqrt{\theta}} \sqrt{\frac{X}{\psi^{-1}(X^{1/4})}},
\]
then partition \(J_{\theta}(X)\) by restricting to terms where \(m, n \leq N\) satisfies \(\left| a_{m,n} \sqrt{X} \right| \leq T\):
\begin{multline*}
J_{\theta}(X) = 
    \frac{X^{3/2}}{2\pi^{2}} 
    \sum_{\substack{
        m,n \leq N \\ 
        |a_{m,n}\sqrt{X}| \leq T
    }} 
    \frac{\tau(m)\tau(n)}{(mn)^{3/4}} \Lambda(a_{m,n}\sqrt{X}) \\
    + \frac{X^{3/2}}{2\pi^{2}} 
    \sum_{\substack{
        m,n \leq N \\ 
        |a_{m,n}\sqrt{X}| > T
    }} 
    \frac{\tau(m)\tau(n)}{(mn)^{3/4}} \Lambda(a_{m,n}\sqrt{X}).
\end{multline*}

    For the error term \(\tilde{R}_N(X)\), we have the following estimate:
\begin{align}\label{equation errormain}
    \tilde{R}_N(X) \ll X^{11/8 + \varepsilon}.
\end{align}

Next, we consider the \textit{upper diagonal} terms, which are given by:
\begin{align}\label{equation udiagonal}
    D^u(X) := \frac{X^{3/2}}{2\pi^2} \sum_{\substack{m, n \leq N \\ |a_{m,n}\sqrt{X}| > T}} \frac{\tau(m)\tau(n)}{(mn)^{3/4}} \Lambda(a_{m,n}\sqrt{X}).
\end{align}
We use \(\left|\Lambda(u)\right| \ll |u|^{-1}\) to bound the terms. Therefore, we obtain the following:
\[
D^u(X) \ll \frac{X^{3/2}}{T} \sum_{m, n \leq N} \frac{\tau(m)\tau(n)}{(mn)^{3/4}}.
\]
This simplifies to
\[
D^u(X) \ll X \cdot \psi^{-1}(X^{1/4})^{1/2} \left( \sum_{m \leq N} \frac{\tau(m)}{m^{3/4}} \right)^2.
\]
For the next step, we define the function \(A(t)\) as
\[
A(t) := \sum_{m \leq t} \tau(m) \quad (t \geq 1).
\]
We now apply the Riemann-Stieltjes integral together with the identity from \eqref{equation hyperbola} to expand the following equality:
\[
\sum_{m \leq N} \frac{\tau(m)}{m^{3/4}} = \int_1^N \frac{1}{t^{3/4}}  dA(t).
\]
Using \(\Delta(x)\ll x\), we can asymptotically bound the above integral as follows:
\[
N^{1/4} \log N + (2\gamma - 1) N^{1/4} + \frac{\Delta(N)}{N^{3/4}} + \frac{3}{4}\int_1^N \left( \frac{\log t}{t^{3/4}} + \frac{2\gamma - 1}{t^{3/4}} + \frac{\Delta(t)}{t^{7/4}} \right)  dt.
\]
This simplifies to
\[
\int_1^N \frac{1}{t^{3/4}}  dA(t) \ll N^{1/4} \log N.
\]
Therefore, we conclude that\[
D^u(X) \ll X^{11/8} \left( \log X \right)^2 \psi^{-1}(X^{1/4})^{1/2}.
\]

Next, for the \textit{lower diagonal} terms, we define:
\begin{align}\label{equation ldiagonal}
   D^l(X) := \frac{X^{3/2}}{2\pi^2} \sum_{\substack{m, n \leq N \\ |a_{m,n}\sqrt{X}| \leq T}} \frac{\tau(m) \tau(n)}{(mn)^{3/4}} \Lambda(a_{m,n} \sqrt{X}).
\end{align}
Since \(\psi^{-1}(X^{1/4}) \geq 2\) for sufficiently large \(X \geq 1\) (it is increasing and unbounded), and \(\left|n - m\theta\right| \leq \frac{1}{2}\sqrt{\frac{m}{\psi^{-1}(X^{1/4})}} + \frac{1}{2}\) whenever \(\left|a_{m,n}\sqrt{X}\right| \leq T\), we define:
\[
U := \frac{1}{2} \sqrt{\frac{m}{\psi^{-1}(X^{1/4})}} + \frac{1}{2}.
\]
Using this definition, we estimate \(D^l(X)\) as follows:
\[
D^l(X) \ll X^{3/2} \sum_{m \leq N} \frac{\tau(m)}{m^{3/4}} \sum_{\substack{n \leq N \\ |n - m\theta| \leq U}} \frac{\tau(n)}{n^{3/4}} \left|\Lambda(a_{m,n} \sqrt{X})\right|.
\]
We now partition the sum according to the value of \(U\). Specifically, we observe that \(U < 1\) if and only if \(m < \psi^{-1}(X^{1/4})\). For the terms where \(U < 1\),  we can write their contribution as:
\[
D^l_{U < 1}(X) := X^{3/2} \sum_{\substack{m < \psi^{-1}(X^{1/4})}} \frac{\tau(m)}{m^{3/4}} \sum_{\substack{n \leq N \\ |n - m\theta| \leq U < 1}} \frac{\tau(n)}{n^{3/4}} \left|\Lambda(a_{m,n} \sqrt{X})\right|.
\]Using \(\left|\Lambda(u)\right| \ll |u|^{-1}\), this can be estimated as:
\[
D^l_{U < 1}(X) \ll X N^\varepsilon \sum_{\substack{m < \psi^{-1}(X^{1/4})}} \frac{1}{m^{3/4}} \sum_{\substack{n \leq N \\ |n - m\theta| \leq U < 1}} \frac{\sqrt{n} + \sqrt{m\theta}}{n^{3/4} \left|n - m\theta\right|}.
\]
Furthermore, employing the Diophantine properties of \(\theta\) and noting that \(U \ll \sqrt{m}\), we observe that:
\[
|n - m\theta| \geq \|m\theta\| \geq \frac{C}{\psi(m)},
\]
and
\[
\frac{m\theta}{2} \leq n \ll m
\]
for all but finitely many \(m\) and \(n\); for the exceptional \(m, n\), the sum above is \(O(X^{1+\varepsilon})\). We notice that for non-exceptional terms we have \(\psi(m) < X^{1/4}\), since \(\psi\) is increasing, which gives:
\[
D^l_{U<1}(X) \ll X N^{\varepsilon} \sum_{m < \psi^{-1}(X^{1/4})} \frac{1}{m} \sum_{\substack{n \leq N \\ |n - m\theta| \leq U < 1}} \psi(m) + X^{1+\varepsilon},
\]
and we conclude that:
\[
D^l_{U<1}(X) \ll X^{5/4 + \varepsilon} \log \left(\psi^{-1}(X^{1/4})\right).
\]
Finally, for the terms on the lower diagonal where \(U \geq 1\), we further partition them into those belonging to the sequence of convergents of \(\theta\) and those that do not:
\[
D^l_{U \geq 1}(X) := X^{3/2} \sum_{\psi^{-1}(X^{1/4}) \leq m \leq N} \frac{\tau(m)}{m^{3/4}} \sum_{\substack{n \leq N \\ |n - m\theta| \leq U}} \frac{\tau(n)}{n^{3/4}} |\Lambda(a_{m,n} \sqrt{X})|.
\]
We split this into two parts:
\begin{align*}
D^l_{U \geq 1}(X) &= X^{3/2} \sum_{\psi^{-1}(X^{1/4}) \leq m \leq N} \frac{\tau(m)}{m^{3/4}} \sum_{\substack{n \leq N \\ |n - m\theta| < \frac{1}{2m}}} \frac{\tau(n)}{n^{3/4}} |\Lambda(a_{m,n} \sqrt{X})| \\
&\quad + X^{3/2} \sum_{\psi^{-1}(X^{1/4}) \leq m \leq N} \frac{\tau(m)}{m^{3/4}} \sum_{\substack{n \leq N \\ \frac{1}{2m} \leq |n - m\theta| \leq U}} \frac{\tau(n)}{n^{3/4}} |\Lambda(a_{m,n} \sqrt{X})|.
\end{align*}

\begin{remark}\label{remark 4.1}
Let \(\theta = [a_0; a_1, a_2, \dots]\) be the simple continued fraction expansion of \(\theta\), and let \(\left( n_k/m_k \right)_{k=0}^\infty\) denote its convergents. The recurrence relation \(m_{k+1} = a_{k+1} m_k + m_{k-1} \) implies by induction that \(m_k \geq F_{k+1}\), where \(F_{k+1}\) is the \((k+1)\)th Fibonacci number. Applying Binet's formula, \(F_k = \frac{\phi^k - (-\phi)^{-k}}{\sqrt{5}}\) with \(\phi = \frac{1 + \sqrt{5}}{2}\), we obtain the exponential lower bound \(m_k \geq \phi^{k-1}\) for all \(k \geq 1\).
\end{remark}

For terms satisfying \(|n - m\theta| < \frac{1}{2m}\) (corresponding to the first summation term), we apply Lemma \ref{lemma Legendre}. This lemma establishes that \(n/m\) must be a convergent of \(\theta\), implying that the denominator \(m\) grows at least exponentially. Using the boundedness of \(\Lambda\) (from Lemma \ref{lemma functionLambda}), \(\tau(m) = o(m^\varepsilon)\) and the exponential lower bound \(m_k \geq \phi^{k-1}\) (from Remark \ref{remark 4.1}), we derive the following estimate:
\[
X^{3/2} \sum_{\psi^{-1}(X^{1/4}) \leq m \leq N} \frac{\tau(m)}{m^{3/4}} \sum_{\substack{n \leq N \\ |n - m\theta| < \frac{1}{2m}}} \frac{\tau(n)}{n^{3/4}} |\Lambda(a_{m,n} \sqrt{X})|
\]
\[
\ll X^{3/2} \sum_{\psi^{-1}(X^{1/4}) \leq m_k \leq N} \frac{1}{m_k^{3/2 - \varepsilon}} \ll X^{3/2} \sum_{\phi^{k-1} \geq \psi^{-1}(X^{1/4})} \phi^{-(3/2-\varepsilon)(k-1)}
\]
\[
\ll \frac{X^{3/2}}{\psi^{-1}(X^{1/4})^{3/2 - \varepsilon}}.
\]
Next, for the remaining terms such that \(\frac{1}{2m} \leq |n - m\theta| \leq U\), corresponding to the last summand in the expression for \(D^l_{U \geq 1}(X)\):
\[
S(X) := X^{3/2} \sum_{\psi^{-1}(X^{1/4}) \leq m \leq N} \frac{\tau(m)}{m^{3/4}} \sum_{\substack{n \leq N \\ \frac{1}{2m} \leq |n - m\theta| \leq U}} \frac{\tau(n)}{n^{3/4}} |\Lambda(a_{m,n} \sqrt{X})|.
\]
This expression splits into two parts:
\begin{align*}
S(X) &= X^{3/2} \sum_{\psi^{-1}(X^{1/4}) \leq m \leq N} \frac{\tau(m)}{m^{3/4}} \sum_{\substack{n \leq N \\ \frac{1}{2m} \leq |n - m\theta| < 1}} \frac{\tau(n)}{n^{3/4}} |\Lambda(a_{m,n} \sqrt{X})| \\
&\quad + X^{3/2} \sum_{\psi^{-1}(X^{1/4}) \leq m \leq N} \frac{\tau(m)}{m^{3/4}} \sum_{\substack{n \leq N \\ 1 \leq |n - m\theta| \leq U}} \frac{\tau(n)}{n^{3/4}} |\Lambda(a_{m,n} \sqrt{X})|.
\end{align*}
For the last summand in \(S(X)\), we use \(U \ll \frac{\sqrt{m}}{\psi^{-1}(X^{1/4})^{1/2}}\) and \(\left|\Lambda(u)\right| \ll |u|^{-1}\), resulting in the following:
\[
X^{3/2} \sum_{\psi^{-1}(X^{1/4}) \leq m \leq N} \frac{\tau(m)}{m^{3/4}} \sum_{\substack{n \leq N \\ 1 \leq |n - m\theta| \leq U}} \frac{\tau(n)}{n^{3/4}} |\Lambda(a_{m,n} \sqrt{X})|
\]
\[
\ll \frac{X N^{\varepsilon}}{\psi^{-1}(X^{1/4})^{1/2}} \sum_{\psi^{-1}(X^{1/4}) \leq m \leq N} \frac{1}{\sqrt{m}} 
\]
\[
\ll \frac{X^{1 + \varepsilon}}{\psi^{-1}(X^{1/4})}.
\]
Now, for the first summand in \(S(X)\), we carefully partition the sum over \(m\) into two cases: \(X^{1/4} \leq m \leq N\) and \(\psi^{-1}(X^{1/4}) \leq m < X^{1/4}\):
\[
X^{3/2} \sum_{X^{1/4} \leq m \leq N} \frac{\tau(m)}{m^{3/4}} \sum_{\substack{n \leq N \\ \frac{1}{2m} \leq |n - m\theta| < 1}} \frac{\tau(n)}{n^{3/4}} |\Lambda(a_{m,n} \sqrt{X})|
\]
\[
\ll X^{3/2} N^{\varepsilon} \sum_{X^{1/4} \leq m \leq N} \frac{1}{m^{3/2}} 
\]
\[
\ll X^{11/8 + \varepsilon}.
\]
In this step, we used the boundedness of \(\Lambda\). Finally, for the terms where \(\psi^{-1}(X^{1/4}) \leq m < X^{1/4}\), we again apply the bound \(\left|\Lambda(u)\right| \ll |u|^{-1}\):
\[
X^{3/2} \sum_{\psi^{-1}(X^{1/4}) \leq m < X^{1/4}} \frac{\tau(m)}{m^{3/4}} \sum_{\substack{n \leq N \\ \frac{1}{2m} \leq |n - m\theta| < 1}} \frac{\tau(n)}{n^{3/4}} |\Lambda(a_{m,n} \sqrt{X})|
\]
\[
\ll X N^{\varepsilon} \sum_{\psi^{-1}(X^{1/4}) \leq m < X^{1/4}} 1
\]
\[
\ll X^{5/4 + \varepsilon}.
\]
By collecting the estimates for \(S(X)\), we obtain:
\[
S(X) \ll \frac{X^{1 + \varepsilon}}{\psi^{-1}(X^{1/4})} + X^{11/8 + \varepsilon}.
\]
Combining the estimates for \(D^l(X)\), we conclude:
\[
D^l_{U \geq 1}(X) \ll \frac{X^{3/2}}{\psi^{-1}(X^{1/4})^{3/2 - \varepsilon}} + X^{11/8 + \varepsilon}
\]
and
\[
D^l_{U < 1}(X) \ll X^{5/4 + \varepsilon} \log \left(\psi^{-1}(X^{1/4})\right).
\]

Thus, the final estimate for the lower diagonal terms in \eqref{equation ldiagonal} is
\[
D^l(X) \ll \frac{X^{3/2}}{\psi^{-1}(X^{1/4})^{3/2 - \varepsilon}} + X^{11/8 + \varepsilon} + X^{5/4 + \varepsilon} \log \left(\psi^{-1}(X^{1/4})\right).
\]
To complete the proof, we collect the bounds for expressions \eqref{equation errormain}, \eqref{equation udiagonal}, and \eqref{equation ldiagonal}, and use the hypothesis that \(\psi^{-1}(X) = O(X^{1/4})\) to establish the following bound for the integral in \eqref{equation mainintegral}:
\[
\int_1^X \Delta(x) \Delta(\theta x)  dx \ll \frac{X^{3/2}}{\psi^{-1}(X^{1/4})^{3/2 - \varepsilon}} + X^{11/8 + \varepsilon} + X^{11/8} (\log X)^2 \psi^{-1}(X^{1/4})^{1/2}
\]
\[
\ll \frac{X^{3/2}}{\psi^{-1}(X^{1/4})^{3/2 - \varepsilon}}.
\]

The proof of Theorem \ref{theorem main} is complete.
\end{proof}

Now, we give a proof of Theorem \ref{theorem second}:
\begin{proof}[Proof of Theorem 1.2]
    Applying the definition of a finite irrationality measure from \cite{Aymone_etal_2024}, there exists \(C > 0\) such that for any fixed \(\delta > 0\), the inequality
\[
\|m\theta\| \geq \frac{C}{m^{\eta + \delta}}
\]
is satisfied for all but finitely many integers \(m \geq 1\). Equivalently, \(\theta\) is not approximable to order \(\psi(x) = x^{\eta + \delta}/C\).

Let \(\varepsilon > 0\). Following the proof of Theorem \ref{theorem main}, the asymptotic behavior of the integral in \eqref{equation mainintegral} is governed by the bounds for \(D^{l}_{U \geq 1}(X)\) and \(D^{u}(X)\). We show that \(D^{l}_{U \geq 1}(X)\) admits a sharper estimate.

Let \(\left( n_k/m_k \right)_{k=0}^\infty\) denote the convergents of \(\theta\). Applying the bounds \(|\Lambda(u)| \ll |u|^{-1}\), \(\sqrt{m} \ll N^{1/2}\), and \(\tau(m) = o(m^{\varepsilon})\), we bound the principal term of \(D^{l}_{U \geq 1}(X)\):
\[
X^{3/2} \sum_{\substack{\psi^{-1}(X^{1/4}) \leq m \leq N}} \frac{\tau(m)}{m^{3/4}} 
\sum_{\substack{n \leq N \\ |n - m\theta| < \frac{1}{2m}}} \frac{\tau(n)}{n^{3/4}} \left|\Lambda\left(a_{m,n} \sqrt{X}\right)\right|
\]
\[
\ll X^{11/8} \sum_{\substack{\psi^{-1}(X^{1/4}) \leq m_k \leq N}} \frac{1}{m_k^{3/2 - \varepsilon} \| m_k \theta \|}.
\]

From the definition of irrationality measure:
\[
\|m\theta\| \geq \frac{C}{m^{\eta + \varepsilon}}
\]
for sufficiently large integers \(m \geq 1\), which implies \(\frac{1}{||m_k\theta||} \ll m_k^{\eta + \varepsilon}\). Combining this with the exponential growth \(m_k \geq \phi^{k-1}\) (Remark \ref{remark 4.1}) yields:
\[
X^{11/8} \sum_{\substack{\psi^{-1}(X^{1/4}) \leq m_k \leq N}} \frac{1}{m_k^{3/2 - \eta - \varepsilon}} 
\ll X^{11/8} \sum_{\substack{\phi^{k-1} \geq \psi^{-1}(X^{1/4})}} \phi^{-(3/2 - \eta - \varepsilon)(k-1)}.
\]
Since \(\eta < 3/2\), the geometric series converges, giving:
\[
\ll X^{11/8}.
\]
Thus we obtain:
\[
D^{l}_{U \geq 1}(X) \ll X^{11/8 + \varepsilon}.
\]

Combining all estimates for \eqref{equation mainintegral} as in Theorem \ref{theorem main}, we select \(\delta > 0\) satisfying \(\frac{1}{8(\eta + \delta)} < \varepsilon\):
\[
\int_1^X \Delta(x) \Delta(\theta x)  dx 
\ll X^{11/8 + \varepsilon} 
+ X^{11/8} (\log X)^2 \psi^{-1}\bigl(X^{1/4}\bigr)^{1/2} 
\ll X^{3/2 - 1/8 + \varepsilon}.
\]
This completes the proof of Theorem \ref{theorem second}.
\end{proof}

\section{Concluding Remarks and Open Questions}
We have established an upper bound for \eqref{equation mainintegral} in terms of \(\psi^{-1}\) for irrational \(\theta > 0\) with the prescribed irrationality measure \(\psi\), including an improved estimate for the generic case. For any Liouville number \(\tau_\beta > 0\) with finite irrationality base \(\beta \geq 1\) (Example \ref{example 1.1}), we obtained the bound
\[
O_{\varepsilon}\left( \frac{X^{3/2}}{(\log X)^{3/2 - \varepsilon}} \right).
\]

Significant open problems remain:

\textbf{Problem 1.} For Liouville numbers with infinite irrationality base \\
(extremely well-approximable), can one establish bounds sharper than \\
\(o(X^{3/2})\) beyond specific constructions like Example \ref{example 1.2}?

\textbf{Problem 2.} What \(\Omega\)-bounds (lower bounds) exist for \eqref{equation mainintegral} \\
across different classes of irrational \(\theta\)?

\noindent\textbf{Acknowledgements.} The author is deeply grateful to his supervisor, Marco Aymone, for proposing this research problem, offering invaluable guidance throughout its development, and for his meticulous reading of earlier drafts of this manuscript. The author also wishes to thank the anonymous referee for his insightful suggestions, which greatly improved the article.

This work was supported by the Research Support Foundation of the State of Minas Gerais (FAPEMIG) under the Graduate Support Program (PAPG) [Master's Scholarship in Mathematics, Grant No. 328].

\bibliographystyle{ieeetr}
\bibliography{main.bib}

\end{document}